\documentclass[11pt]{amsart}
\usepackage{amsmath,amssymb,amsfonts,epsf,epsfig}
\usepackage{enumerate,amsthm,color,mathdots}
\usepackage{euscript,graphicx}
\usepackage{url, float, hyperref, color} 
\usepackage[usenames,dvipsnames,svgnames,table]{xcolor}
\usepackage{graphicx, mathtools}
\usepackage{tikz}
\usepackage{lineno}

\newtheorem{theorem}{Theorem}[section]
\newtheorem{lemma}[theorem]{Lemma}
\newtheorem{proposition}[theorem]{Proposition}
\newtheorem{corollary}[theorem]{Corollary}

\newtheorem{definition}[theorem]{Definition}

\newtheorem{remark}[theorem]{Remark}

\renewcommand{\mod}{\hbox{{\rm mod}}\, }
\newcommand{\Aut}{{\rm Aut}}

\newcommand{\Sym}{{\rm S}}
\newcommand{\psl}{{\rm PSL}}
\newcommand{\pgl}{{\rm PGL}}
\newcommand{\Sz}{{\rm Sz}}

\newcommand{\ZZ}{\mathbb{Z}}

\newcommand{\F}{{\mathcal{F}}}
\newcommand{\M}{{\mathcal{M}}}
\newcommand{\Mon}{\hbox{\rm Mon}}

\newcommand{\rD}{\hbox{\rm D}}
\newcommand{\rS}{\hbox{\rm S}}
\newcommand{\rA}{\hbox{\rm A}}

\title[Projective maniplexes]{Maniplexes with automorphism group $\textrm{PSL}(2,q)$}
\author{Dimitri Leemans}
\email{leemans.dimitri@ulb.be}
\author{Micael Toledo}
\email{micaelalexitoledo@gmail.com}
\address{Universit\'e Libre de Bruxelles, D\'epartement de Math\'ematique, C.P.216 - Alg\`ebre et Combinatoire, Boulevard du Triomphe, 1050 Brussels, Belgium}
\thanks{The authors gratefully acknowledge financial support from the F\'ed\'eration Wallonie-Bruxelles -- Actions de Recherche Concert\'ees (ARC Advanced grant)}

\begin{document}

\begin{abstract}
A maniplex of rank $n$ is a combinatorial object that generalises the notion of a rank $n$ abstract polytope. A maniplex with the highest possible degree of symmetry is called reflexible. In this paper we prove that there is a rank $4$ reflexible maniplex with automorphism group $\psl_2(q)$ for infinitely many prime powers $q$, and that no reflexible maniplex of rank $n > 4$ exists that has $\psl_2(q)$ as its full automorphism group. 
\end{abstract}
\maketitle


\section{Introduction}

A maniplex is a combinatorial object that generalises the notion of a map on a surface to higher dimensions. Thus, maniplexes also generalise abstract polytopes in a similar way that maps generalise abstract and geometric polyhedra. As with maps and abstract polytopes, the study of maniplexes is essentially that of their symmetries, with much of the research focusing on those objects exhibiting the highest possible degree of symmetry. These are typically called regular in the case of abstract polytopes, and reflexible, in the case of maps and maniplexes. Informally, an abstract regular polytope has an automorphism group that is called a string C-group, namely a string group generated by involutions satisfying an intersection condition that ensures, among others, that the polytope is strongly connected. Meanwhile, the automorphism group of a reflexible maniplex of rank $n$ is also a string group generated by involutions, but we do no require it to satisfy the intersection property. In this sense every regular abstract polytope is a reflexible maniplex. However, this generalisation extends beyond regularity or reflexibility: every abstract polytope, regardless of its symmetry, is a maniplex. A charaterization of polytopality for maniplexes was given in \cite{polymani}. Maniplexes were first introduced in \cite{manisteve}, and have since received considerable attention (see for instance \cite{stg,mani2,2orbit,dansteve}). 


 Maps and abstract polytopes having a projective special linear group as its automorphism group have been extensively studied. In particular, Sjerve and Cherkassov showed in \cite{cherk}  that if $q \neq 2,3,7,9$, then there exists a reflexible map with automorphism group $\psl_2(q)$. A complete enumeration of maps (in the wider context of regular hypermaps) having $\psl_2(q)$ as a subgroup of index $k \leq 2$ of their automorphism group was given in \cite{Conder2008} by Conder et al., extending the work of Macbeath \cite{macbeath} and Sah \cite{sah}. Roughly at the same time, Leemans and Schulte \cite{pslpoly} showed that there exist exactly one regular abstract polytopes of rank $4$ with automorphism group $\psl_2(q)$ for each $q \in \{11,19\}$, and that these are the only two regular abstract polytopes of rank $4$ with a projective special linear group as their full automorphism group. Moreover, Leemans and Vauthier proved in~\cite{LV2005} that for $k \geq 5$, there are no regular abstract polytopes of rank $k$ with automorphism group $\psl_2(q)$ for any primer power $q$.
 
 In this paper, we extend the results of \cite{cherk} and \cite{pslpoly} by proving that there exists a reflexible maniplex of rank $4$ with automorphism group $\psl_2(q)$ for infinitely many odd primes $q$ and that no such maniplex exists if $q$ is an even prime power. We show that any such maniplex is either isomorphic to one of the two (polytopal) maniplexes of \cite{pslpoly}, or it has a single facet (alt. vertex). In the latter case, one of the two groups associated with the facet and the vertex-figure is isomorphic to the whole automorphism group $\psl_2(q)$ while the other is isomorphic to either $\rA_5$ or $\psl_2(q)$. Therefore such a maniplex (or its dual) has one facet and either one or $|\psl_2(q)|/60$ vertices. We construct an infinite family of maniplexes for each of the two aforementioned possibilities. We also show that there are no reflexible maniplexes of rank $k \geq 5$ with automorphim group $\psl_2(q)$

 In Section \ref{sec:def} we formally define maniplexes and their relation to string groups generated by involutions. In Section \ref{sec:string} we prove several auxiliary results on the possible string group representation for the groups $\psl_2(q)$. In Section \ref{sec:main} we construct two infinite families of rank $4$ maniplexes with automorphism group $\psl_2(q)$ where $q$ is an odd prime, and we show that no such maniplex exist if $q$ is an even prime power. 
 Finally, in Section~\ref{sec:remarks} we give some concluding remarks.






\section{Preliminary concepts}
\label{sec:def}

Given a set $\Omega$, we denote by $\Sym_{\Omega}$ the group of all permutations of the set $\Omega$. We denote by $id$ the permutation of $\Sym_\Omega$ fixing all elements of $\Omega$. Sometimes we also write $\Sym_n$ where $n = |\Omega|$ when the context is clear.
For a prime power $q$, we let $\psl_2(q)$, $\pgl_2(q)$ and $\textrm{SL}_2(q)$ denote the projective special linear group, the projective general linear group and the special linear group of the vector space of dimension $2$ over the finite field of order $q$, respectively. We let $\rD_{2k}$ denote the dihedral group of order $2k$.

For a positive integer $n$, a {\em maniplex of rank} $n$ (also called an {\em $n$-maniplex}) is a pair $\M := (\F, (r_0,r_1,\ldots,r_{n-1}))$ where $\F$ is a non-empty set, whose elements are called {\em flags}, and each $r_i$ is an involution in $\Sym_{\F}$ with the following properties:
\begin{enumerate}
\item $r_i$ and $r_ir_j$ are fixed-point-free for all $i$ and $j \neq i$;
\item $(r_ir_j)^2 = id$ whenever $|i-j| > 1$;
\item the action of the group $\Mon(\M)= \langle r_0, r_1, \ldots, r_{n-1} \rangle$ is transitive on $\F$.
\end{enumerate} 

The group $\Mon(\M)$ is called the {\em monodromy group} of $\M$ and the $n$-tuple $(r_0,\ldots,r_{n-1})$ is the (ordered) set of {\em distinguished generators} of $\Mon(\M)$. For a flag $\Phi$ and a permutation $r \in \Mon(\M)$, we let $\Phi^r$ denote the image of $\Phi$ under $r$. In particular, for $i \in \{0,\ldots,n-1\}$ we call the flag $\Phi^{r_i}$ the $i$-neighbour of $\Phi$. Since each $r_i$ is an involution, the ``neighbour" relation is symmetric. If two flags $\Phi$ and $\Psi$ are $i$-neighbours of one another we also say that they are {\em $i$-adjacent}.

\begin{remark}
A rank $3$ maniplex coincides with the notion of a map on a surface, and maps are sometimes defined in purely combinational terms as a set of three involutory permutations acting on a set (see \cite{combimaps} for example, where flags are called `blades'). Given that the notion of a map is  fairly well-established, we will refer to $3$-maniplexes simply as maps.   
\end{remark}

For $i \in \{0,\ldots,n-1\}$ an {\em $i$-face} of $\M$ is an orbit of the subgroup $G_i = \langle r_j \mid j \neq i \rangle$ of $\Mon(\M)$ on its action on $\F$. We call $0$- and $(n-1)$-face a {\em vertex} and a {\em facet}, respectively. 

Consider the group $\Mon(\M)^+ = \langle r_ir_j \mid i,j \in \{0,\ldots,n-1\} \rangle$ and note that it has at most two orbits on flags. We say $\M$ is {\em non-orientable} if $\Mon(\M)^+$ has a single orbit, and we say it is {\em orientable} otherwise. In the latter case, each of the two orbits is called an {\em orientation}.

An automorphism of $\M$ is a permutation of the flags that maps $i$-adjacent flags to $i$-adjacent flags for all $i \in \{0,\ldots,n-1\}$. The group of automorphisms of $\M$ is denoted $\Aut(\M)$. Since automorphisms must preserve $i$-adjacencies, the actions of $\Aut(\M)$ and $\Mon(\M)$ commute. That is, for every $\Phi \in \F$, $\alpha \in \Aut(\M)$ and $r \in \Mon(\M)$ we have

\begin{align}
\Phi^{\alpha r} =  \Phi^{r \alpha}.
\end{align}

From this fact and from the transitivity of $\Mon(\M)$, we see that an automorphism of $\M$ is completely determined by its image on a single flag and thus the action of $\Aut(\M)$ is semiregular on $\F$. We say $\M$ is {\em reflexible} if $\Aut(\M)$ acts transitively on the flags of $\M$. Note that in this case the action of $\Aut(\M)$ on the flags of $\M$ must also be regular, and thus for any two flags there exists a unique automorphism mapping one to the other.

\subsection{String representations}

We say a group $G$ is a {\em string group generated by involutions}, or a $sggi$ for short, if $G$ is generated by involutions $g_0, \ldots, g_{n-1}$ such that $g_ig_j = g_jg_i$ if $|i-j| > 1$. If there are $n$ pairwise distinct such involutions, then we say $G$ is of rank $n$ and the $n$-tuple $(g_0,\ldots,g_{n-1})$ is called a {\em string representation} of $G$ of rank $n$. Two string representations $(g_0,\ldots,g_{n-1})$ and $(g_0',\ldots,g_{n-1}')$ of $G$ are said to be isomorphic if there exists a group automorphism $\phi \colon G \to G$ mapping each $g_i$ to $g_i'$.

As it transpires, the automorphism group of a reflexible $n$-maniplex is a $sggi$ of rank $n$, and every rank $n$ $sggi$ is the automorphism group of a reflexible $n$-maniplex. Then, determining the existence of rank $n$-maniplexes with a given group of automorphisms $G$ reduces to determining the possible rank $n$ string representations of $G$. Let us discuss the relation between these two notions a little bit further.

Let $\M$  be a reflexible $n$ maniplex, let $I = \{0,\ldots,n-1\}$ and let  $\Phi$ be a fixed flag of $\M$. For $i \in I$ let $\rho_i$ be the unique automorphism of $\M$ such that $\Phi^{\rho_i} = \Phi^i$. If $\alpha \in \Aut(\M)$ then, by the transitivity of $\Mon(\M)$, there exist $i_0, i_1,\ldots,i_j \in I$ such that $\Phi^\alpha = \Phi^{r_{i_0}r_{i_1}\ldots r_{i_j}}$. Since the action of $\Aut(\M)$ and $\Mon(\M)$ commute we have

\begin{align*}
\Phi^\alpha = \Phi^{r_{i_0}r_{i_1}\ldots r_{i_j}} = (\Phi^{\rho_{i_0}})^{r_{i_1}\ldots r_{i_j}} = (\Phi^{r_{i_1}\ldots r_{i_j}})^{\rho_{i_0}}=\ldots = \Phi^{\rho_{i_j}\rho_{i_{j-1}}\ldots\rho_{i_0}}, 
\end{align*}

and thus $\alpha = \rho_{i_j}\rho_{i_{j-1}}\ldots\rho_{i_0}$. It follows that $\{\rho_i \mid i \in I\}$ is a generating set for $\Aut(\M)$. Moreover, if $|i-j| > 1$, then $\Phi^{\rho_i\rho_j} = \Phi^{r_jr_i} = \Phi^{r_ir_j} = \Phi^{\rho_j\rho_i}$ and thus $\rho_i\rho_j = \rho_j\rho_i$. Then the $n$-tuple $(\rho_0,\rho_1,\ldots,\rho_{n-1})$ is a string representation of $\Aut(\M)$. Note that we have some freedom in choosing the base flag $\Phi$, and thus this string representation is not unique. However, all the different string representations obtained by considering different flags of $\M$ are conjugate in $G$ and are therefore isomorphic.

Conversely let $G$ be a group with a string representation $(g_0,g_1,\ldots,g_{n-1})$. For $i \in \{0,\ldots,n-1\}$, define $r_i$ as the permutation induced by $g_i$ on $G$ by left multiplication. It is routine to check that $M = (G,(r_0,r_1,\ldots,r_{n-1}))$ is a regular $n$-maniplex with $ \Aut(\M) \cong \Mon(\M) \cong G$, and the actions of $\Mon(\M)$ and $\Aut(\M)$ are isomorphic to the left and right regular actions of $G$, respectively. Furthermore, if we let $\Phi = 1_G$ then the string representation of $G$ associated to $\M$ and $\Phi$ is precisely $(g_0,\ldots,g_{n-1})$.  

By the discussion above, we may slightly abuse language and redefine reflexible $n$-maniplexes in terms of string representations.

\begin{definition}
A reflexible $n$-maniplex $\M$ is a pair $(G,(\rho_0,\rho_1,\ldots,\rho_{n-1}\})$ where $G$ is a group and $(\rho_0,\rho_1,\ldots,\rho_{n-1})$ is a string representation of $G$.
\end{definition}

Let $\sigma_i$ denote the order of $r_ir_{i+1}$ for $i \in \{0,\ldots,n-2\}$. The {\em type} of $\M$ is $\{\sigma_0,\sigma_1,\ldots,\sigma_{n-2}\}$.

\section{String representations of $\psl_2(q)$}
\label{sec:string}

The determination of the subgroup lattice of $\psl_2(q)$ due to Dickson \cite{dickson} (and revisited by many authors such as~\cite{Hup67, Van85}) has become a central tool in studying the possible string representations of projective special linear groups. In \cite{cherk}, it was shown that $\psl_2(q)$ admits a rank $3$ string representation whenever $q \notin \{2,3,7,9\}$. Later Leemans and Schulte \cite{pslpoly} proved that if $q \neq 11,19$, then $\psl_2(q)$ does not admit a rank $4$ string representation $(\rho_0,\rho_1,\rho_2,\rho_3)$ satisfying the following additional {\em intersection property}:

\begin{align}
\label{eq:ip}
\forall \enspace J,K \subseteq \{0,1,2,3\},\quad \langle \rho_j \mid j \in J \rangle \cap \langle \rho_k \mid k \in K \rangle = \langle \rho_j \mid j \in J \cap K \rangle.
\end{align} 

String groups generated by involutions satisfying the intersection property above are usually called {\em string $C$-groups}, and are associated to regular abstract polytopes in much the same way as an $sggi$ is related to a reflexible maniplex. The intersection property is very restrictive in the sense that not many $sggis$ satisfy it. While $\psl_2(q)$ admits only two rank $4$ string representations with the intersection property (one when $q =11$ and one when $q =19$), it admits a string representation for infinitely many values of $q$.

Note that all such representation are of rank $n \leq 4$, as $\psl_2(q)$ does not contain subgroups that are the direct product of two dihedral subgroups or order $6$ or larger. Indeed, if $n>4$ and $(\rho_0,\ldots,\rho_{n-1})$ is a string representation of $\psl_2(q)$, then $\langle \rho_0, \rho_1 \rangle$ and $\langle \rho_{n-2}, \rho_{n-1} \rangle$ are dihedral of order at least $6$ and they centralize each other, making their direct product a subgroup of $\psl_2(q)$. This was first remarked in~\cite[Theorem 2]{LV2005} for string $C$-group representations of $\psl_2(q)$, but the same argument works here. 


Our strategy to produce rank $4$ string representations for $\psl_2(q)$ will be similar to that used in \cite{cherk} and \cite{pslpoly}. We will make extensive use of the results of Dickson on the subgroup lattice of $\psl_2(q)$, summarized in Lemma \ref{lem:dick} below. More detailed versions of this theorem can be found in \cite{michael,pslpoly, Van85}.
\begin{lemma}
\cite{dickson}
\label{lem:dick}
Let $q = p^r$ where $p$ is a prime power and $r$ is a positive integer, and let $H$ be a proper non-abelian subgroup of $\psl_2(q)$. Then one the following holds:
\begin{enumerate}
\item $H$ is a semidirect product of an elementary abelian group with a cyclic group;
\item $H$ is a dihedral group;
\item $H$ is a projective linear group isomorphic to $\psl_2(p^a)$ with $a \mid r$ or $\pgl_2(p^b)$ with $2b \mid r$;
\item $H$ is isomorphic to $\rA_4$, $\rS_4$ or $\rA_5$.
\end{enumerate}
\end{lemma}

For the remainder of this section let $G \cong \psl_2(q)$, $q >3$, and let $(\rho_0,\rho_1,\rho_2,\rho_3)$ be a string representation of $G$. For each $i \in \{0,1,2,3\}$, let $G_i = \langle \rho_j \mid j \neq i \rangle$. We will often refer to $G_0$ and $G_3$ as the {\em vertex group} and {\em facet group} of $\M$ respectively. Observe that the group $G_3$ cannot be abelian, for otherwise $\rho_0$ would lie in the centre of $G$, which is trivial. Then, if $G_3$ is a proper subgroup of $G$, it must belong to one of the four classes of groups described in Lemma \ref{lem:dick}. The main goal of this section is to show that if $G_3$ is a proper subgroup then it is isomorphic to $\rS_4$ or $\rA_5$, by ruling out, one by one, all the other possibilities. This is essentially the same strategy as that used in \cite{pslpoly}, but in the more general setting of maniplexes (that is, of $sggis$ that do not necessarily satisfy the intersection property). In fact, we show that some of the results of \cite{pslpoly} can be restated to not require the intersection property without changing the original proof too much, provided that we are careful. We will need the following three auxiliary results, which can be traced back in one form or another to \cite{dickson}. Lemma~\ref{lem:centraliser} is obtained by inspecting the list of subgroups of $\psl_2(q)$. 

\begin{lemma} \cite[Lemma 3]{pslpoly}
\label{lem:centre}
The centre of a non-abelian subgroup of $\psl_2(q)$ has order at most $2$.
\end{lemma}

\begin{lemma}
\label{lem:centraliser}
If $G \cong \psl_2(q)$ and $r \in G$ is an involution, then the centralizer $C_G(r)$ is a maximal subgroup of $G$ and is isomorphic to the dihedral group $\rD_{2n}$ for some $n \in \{\frac{q + 1}{2}, \frac{q - 1}{2}\}$.
\end{lemma}

\begin{lemma} \cite[Lemma 6]{pslpoly}
\label{lem:dih}
Let $q = p^s$ for some prime $p$ and $s\geq 1$ and let $G \cong \psl_2(q)$. Let $q' \mid q$ and $n \in \{\frac{q' + 1}{2},\frac{q' - 1}{2}\}$. If $k \mid n$ then the following holds:
\begin{enumerate}
\item every dihedral sugbroup $\rD_{2k} \leq G$ of order $2k$ is contained in exactly one dihedral subgroup $\rD_{2n}$ of order $2n$;
\item every dihedral subgroup of order $2n$ is contained in exactly one subgroup of $G$ isomorphic to $\psl_2(q')$. 
\end{enumerate}
\end{lemma}



We will start by proving that $G_3$ is not the semidirect product of an elementary abelian group with a cyclic group.

\begin{lemma}
If $(\rho_0,\rho_1,\rho_2,\rho_3)$ is a string representation of $\psl_2(q)$ then the order of $\rho_i\rho_{i+1}$ is greater than $2$ for all $i \in \{0,1,2\}$.
\end{lemma}

\begin{proof}
Suppose $\rho_0\rho_1$ has order two. Then $\rho_0$ commutes with $\rho_1$, but also with $\rho_2$ and $\rho_3$, and thus $\rho_0$ lies in the center of $\psl_2(q)$, a contradiction. By the same token, the order of $\rho_2\rho_3$ is greater than two. Finally if $(\rho_1\rho_2)^2 = 1_G$ then $\psl_2(q) = \langle \rho_0,\rho_1 \rangle \times \langle \rho_2,\rho_3 \rangle$, which is impossible since $\psl_2(q)$
is not the direct product of two dihedral subgroups.
\end{proof}

This shows that $G_3$ contains a dihedral subgroup of order at least $6$, namely, $\langle \rho_0,\rho_1 \rangle$. Since a semidirect product of an elementary abelian group with a cyclic group does not contain dihedral groups of order larger than $4$, we have the following corollary.

\begin{corollary}
\label{cor:nosd}
$G_3$ is not a semidirect product of an elementary abelian group with a cyclic group.
\end{corollary}

The next step is to show that $G_3$ cannot be a projective linear group. Compare Lemmas \ref{lem:nopsl} and \ref{lem:nopgl} with Lemmas 7 and 9 of \cite{pslpoly}.

\begin{lemma}
\label{lem:nopsl}
Let $G \cong \psl_2(q)$ and let $(\rho_0,\rho_1,\rho_2,\rho_3)$ be a rank $4$ string representation of $G$. The subgroups $G_0 = \langle \rho_1,\rho_2,\rho_3 \rangle$ and $G_3 = \langle \rho_0,\rho_1,\rho_2 \rangle$ cannot be isomorphic to $\psl_2(q')$ with $q' < q$.
\end{lemma}

\begin{proof}
Suppose for a contradiction that $G_3 \cong \psl_2(q')$ for some $q' < q$, and consider the subgroup $(G_3)^{\rho_3}$. Note that $G_3 \neq (G_3)^{\rho_3}$ by the simplicity of $G$ (as $q$ must be at least $4$). But then the dihedral group $\langle \rho_0, \rho_1 \rangle$ lies in the intersection of $G_3$ and $(G_3)^{\rho_3}$, both of which are isomorphic to $\psl_2(q')$ and, by Lemma~\ref{lem:dih}, $G_3 = (G_3)^{\rho_3}$, a contradiction.
\end{proof}

\begin{lemma}
\label{lem:nopgl}
Let $G \cong \psl_2(q)$ and let $(\rho_0,\rho_1,\rho_2,\rho_3)$ be a rank $4$ string representation of $G$. The subgroups $G_0 = \langle \rho_1,\rho_2,\rho_3 \rangle$ and $G_3 = \langle \rho_0,\rho_1,\rho_2 \rangle$ cannot be isomorphic to $PGL_2(q')$ with $q' < q$.
\end{lemma}

\begin{proof}
Lemma 9 of \cite{pslpoly} states that if $G \cong \psl_2(q)$ and $(\rho_0,\rho_1,\rho_2,\rho_3)$ is a rank $4$ string representation of $G$ satisfying the intersection property (\ref{eq:ip}), then the subgroup $G_3 = \langle \rho_0,\rho_1,\rho_2 \rangle$ cannot be isomorphic to $PGL_2(q')$ with $q' < q$. However, the requirement that $(\rho_0,\rho_1,\rho_2,\rho_3)$ satisfies the intersection property is used only to show that the order of the element $\rho_0\rho_1$ is odd. Therefore, it suffices to show that $\rho_0\rho_1$ has odd order and the result will follow from the proof of \cite[Lemma 9]{pslpoly}, which we omit here for brevity. If $G_3 = G$ then clearly $G_3 \neq PGL_2(q')$ with $q' < q$. Therefore we can assume $G_3$ is a proper subgroup of $G$, and in particular $\rho_3 \notin G_3$. Then the subgroup $\langle \rho_0,\rho_1,\rho_3 \rangle$ is isomorphic to $\rD_{2t} \times \ZZ_2$ where $t$ is the order of $\rho_0\rho_1$. If $t$ is even, then the group $\rD_{2t} \times \ZZ_2$ is a non-abelian subgroup whose centre has more than $2$ elements, contradicting Lemma \ref{lem:centre}. 
\end{proof}

\begin{lemma}
\label{lem:nodih}
Let $G \cong \psl_2(q)$ and let $H \leq G$ be a dihedral subgroup admitting a rank $3$ string representation $(\rho_0,\rho_1,\rho_2)$. If $\rho_3 \in G$ is an involution commuting  with $\rho_0$ and $\rho_1$, then the group $\langle \rho_0,\rho_1,\rho_2,\rho_3 \rangle$ is a proper subgroup of $G$.
\end{lemma}

We give two distinct proofs of this lemma, one based on geometric arguments, the other based on algebraic arguments.
\begin{proof}[Geometric proof]
The group $\psl_2(q)$ acting on its natural representation on $q+1$ points sits inside $\pgl_2(q)$ which is sharply 3-transitive on these $q+1$ points.

Suppose that $q\equiv 1$ $\mod 4$. Then all involutions of $\psl_2(q)$ have two fixed points.
Also, the only maximal dihedral subgroups containing pairs of involutions that commute are $\rD_{q-1}$.
So $H$ must be in a subgroup $\rD_{q-1}$.
The central involution of $\rD_{q-1}$ fixes two points $a$ and $b$.
Every non-central involution swaps these two points and fixes two points of the $q-1$ remaining points.
Say $\rho_i$ fixes $a_i$ and $b_i$.
For $\rho_3$ to commute with $\rho_0$ and $\rho_1$ it needs to swap $a_0$ with $b_0$ and swap $a_1$ with $b_1$.
But the central involution of $\rD_{q-1}$ does that. As $\pgl_2(q)$ is sharply 3-transitive, there is a unique element mapping three points on three other points. Here we have even more, four points. Hence,  $\rho_3$ has to be that central involution and the group $\langle H, \rho_3 \rangle$ is necessarily a subgroup of $\rD_{q-1}$.

Suppose $q\equiv 3$ $\mod 4$. Then $\psl_2(q)$ can be seem as a subgroup of $\psl_2(q^2)$ with $q^2 \equiv 1$ $\mod 4$ and we can apply the same argument as above on $q^2+1$ points.

Finally, suppose that $q$ is even. Then commuting involutions necessarily have a fixed point in common and $H$ must be an elementary abelian group of order $4$. Moreover,   $\langle H, \rho_3 \rangle$ will also necessarily be an elementary abelian group, not $\psl_2(q)$.
\end{proof}

\begin{proof}[Algebraic proof]
Suppose that $\rho_3 \in \psl_2(q)$ is an involution satisfying $(\rho_0\rho_3)^2=(\rho_1\rho_3)^2=1_G$. We will show that $K = \langle \rho_0,\rho_1,\rho_2,\rho_3 \rangle$ cannot be the whole $\psl_2(q)$.
Note that if $\rho_0$ and $\rho_1$ commute, then $\rho_0$ lies in the center of $K$ and so $K \neq \psl_2(q)$. Suppose $\rho_0$ and $\rho_1$ do not commute. Then $D := \langle \rho_0, \rho_1 \rangle$ is a dihedral subgroup of $H = \langle \rho_0,\rho_1,\rho_2\rangle$ of order larger than $4$. Since $\rho_0$ and $\rho_1$ lie in the centraliser $C$ of $\rho_3$ in $\psl_2(q)$, then $D \leq C$. By Lemma \ref{lem:centraliser}, $C$ is a maximal dihedral subgroup of $\psl_2(q)$. Furthermore, by Lemma \ref{lem:dih}, every dihedral subgroup of $\psl_2(q)$ lies inside a unique maximal dihedral subgroup. Let $\bar{H}$ be the maximal dihedral subgroup containing $H$. Then $D \leq \bar{H} \cap C$ and thus $\bar{H} = C$. In particular, this shows that $\rho_2 \in C$ and thus $K = \langle \rho_0,\rho_1,\rho_2,\rho_3 \rangle \leq C \lneq \psl_2(q)$. 
\end{proof}

We now have all the results necessary to prove the main theorem of this section.


\begin{theorem}
\label{the:facet}
Let $q$ be a prime power and let $G \cong \psl_2(q)$. If $(\rho_0,\rho_1,\rho_2,\rho_3)$ is a string representation of $G$, then the subgroup $G_3 = \langle \rho_0, \rho_1, \rho_2 \rangle$ is isomorphic to $\rS_4$, $\rA_5$ or $G$. 
\end{theorem}

\begin{proof}
Suppose $G_3$ is a proper subgroup of $G$. Since $G_3$ is non-abelian, one of the four items of Lemma \ref{lem:dick} must hold for $G_3$. By Corollary \ref{cor:nosd} and Lemmas \ref{lem:nopsl}, \ref{lem:nopgl} and \ref{lem:nodih} we see that $G_3$ must be isomorphic to either $\rA_4$, $\rS_4$ or $\rA_5$. However $\rA_4$ admits no string representation of rank $3$ \cite{cherk}. This leaves only $\rS_4$ and $\rA_5$ as possibilities for $G_3$. 
\end{proof}

\section{Reflexible 4-maniplexes with group $\psl_2(q)$}
\label{sec:main}

Let $\M:=(G,(\rho_0,\rho_1,\rho_2,\rho_3))$ be a reflexible $4$-maniplex with $G \cong \psl_2(q)$ for some prime power $q$. The {\em dual} of $\M$ is the maniplex $\M^* = (G,(\rho_3,\rho_2,\rho_1,\rho_0))$. By applying Theorem \ref{the:facet} to $\M^*$ we see that the group $G_0 = \langle \rho_3,\rho_2,\rho_1 \rangle$ is isomorphic to $\rS_4$, $\rA_5$ or $G$. Hence there are essentially three distinct possibilities for $\M$, depending on what the groups $G_0$ and $G_3$ may look like.

\begin{remark}
\label{cor:classes}
One of the following holds for $\M$:

\begin{enumerate}
\item $G_0 \cong G_3 \cong G$;
\item Exactly one of $G_0$ or $G_3$ is isomorphic to $G$ and the other is isomorphic to either $\rS_4$ or $\rA_5$;
\item Both $G_0$ and $G_3$ are isomorphic to a group in $\{\rS_4,\rA_5\}$.
\end{enumerate}
\end{remark} 
 
If for some $i \in \{1,2,3\}$ item $(i)$ of Remark \ref{cor:classes} above holds for $\M$, then we will say $\M$ belongs to Class $i$. Note that if $q > 3$, then $\psl_2(q)$ is a simple group, and thus $\M$ is necessarily non-orientable, regardless of the class to which it belongs, as its automorphism group contains no subgroup of index $2$.

We will show in Lemma \ref{lem:class1} and Proposition \ref{prop:class2} that Classes (1) and (2) are infinite. Meanwhile, Class (3), which is comprised of exactly two maniplexes, was dealt with in \cite{pslpoly}, in the context of abstract regular polytopes. In particular, \cite[Theorem 3]{pslpoly} states that if $G \cong \psl_2(q)$ and $(\rho_0,\rho_1,\rho_2,\rho_3)$ is a string C-group representation of $G$ such that each of $G_3$ and $G_0$ is isomorphic to either $\rS_4$ or $\rA_5$, then $q \in \{11,19\}$ and for each of the two possible values of $q$, the string C-group representation of $G$ is unique (up to isomorphism). However one can readily check that the intersection property is not used in the proof of \cite[Theorem 3]{pslpoly} and thus the result holds in the more general case where $(\rho_0,\rho_1,\rho_2,\rho_3)$ is a string representation. We then obtain the following proposition (see~\cite{Gr1977,Cox1982} for a definition of the $11$- and $57$-cell respectively).

\begin{proposition}
\label{pro:class3}
If $\M$ is a reflexible $4$-maniplex in Class (3), then $q = 11$ and $\M$ is isomorphic to the $11$-cell or $q =19$ and $\M$ is isomorphic to the $57$-cell. In both cases $G_0 \cong G_3 \cong \rA_5$.
\end{proposition}

Proposition \ref{pro:class3} and Remark \ref{cor:classes} combined give us the following important result.



\begin{lemma}
\label{lem:qeven}
There are no reflexible $4$-maniplexes with automorphism group $\psl_2(2^k)$.
\end{lemma}

\begin{proof}
Let $G = \psl_2(q)$ with $q=2^k$. Note that if $k=1$, then $G = \psl_2(2) \cong \rS_3$ and clearly $G$ does not admit a rank $4$ string representation. We can henceforth assume that $k > 1$ and thus $q \geq 4$, forcing $G$ to be a simple group.

Suppose for a contradiction that $G$ admits a string representation $(\rho_0,\rho_1,\rho_2,\rho_3)$. Recall that when $q$ is even, $\psl_2(q) = \pgl_2(q)$. Since $q$ is even, then by Remark \ref{cor:classes} and Proposition \ref{pro:class3}, we can assume without loss of generality that $G_3 := \langle \rho_0,\rho_1,\rho_2 \rangle$ is isomorphic to $G$.

Consider the natural action of $G$ on $q + 1$ points. Since $q$ is even and $G$ is sharply 3-transitive, every involution in $G$ fixes exactly one point (as if it fixes more points, it must fix at least three and the only element fixing three points is $1_G$), and two involutions that commute must have the same fixed point. Let $x$ be the point fixed by $\rho_0$ and $\rho_3$. As $\rho_0$ (resp. $\rho_1$) and $\rho_2$ (resp. $\rho_3$) commute, $\rho_2$ (resp. $\rho_1$) also fixes $x$. Hence $\langle \rho_0,\rho_1,\rho_2,\rho_3\rangle \leq G_x < G$.
\end{proof}

We say a reflexible $n$-maniplex $\M$ is {\em extendible} if there is an $(n+1)$-maniplex $\M'$ whose facets are all isomorphic to $\M$. The maniplex $\M'$ is then called an {\em extension} of $\M$. 

Let $\M = (G,(\rho_0,\rho_1,\rho_2,\rho_3)$ be a reflexible $4$-maniplex with $G \cong  \psl_2(q)$. Observe that $\M_3 := (G_3,(\rho_0,\rho_1,\rho_2))$ is a reflexible map (a $3$-maniplex) in its own right, and every facet of $\M$ is isomorphic to $\M_3$. Then, $\M$ is an extension of $\M_3$. In particular, in light of Theorem \ref{the:facet}, $\M$ must be an extension of a reflexible map whose full automorphism group is either $\rS_4$, $\rA_5$ or $\psl_2(q)$.

\begin{theorem}
\label{the:ext}
Let $q=p^k$ where $p$ is an odd prime, let $\alpha \in \{-1,1\}$ be such that $q \equiv \alpha$ $(\mod 4)$ and let $n = \frac{q-\alpha}{2}$. For an odd divisor $k$ of $n$, every reflexible map on $\psl_2(q)$ of type $\{k,x\}$, $x \geq 3$, is extendible.
\end{theorem}

\begin{proof}
Let $\M = (\psl_2(q),(\rho_0,\rho_1,\rho_2))$ be a regular map of type $\{k,x\}$ for an odd divisor $k$ of $n$ and some integer  $x\geq 3$. Then the order of $\rho_0\rho_1$ is $k$ and $H := \langle \rho_0, \rho_1 \rangle$ is a dihedral group of order $2k$. By Lemma \ref{lem:dih} $H$ is contained in a subgroup $D < \psl_2(q)$ isomorphic to the dihedral group $\rD_{2n}$. Since $n$ is even, $D$ contains a central involution $\rho_3$, which commutes with $\rho_0$ and $\rho_1$. Note that $\rho_3 \neq \rho_2$, since $\rho_2$ does not commute with $\rho_1$. Moreover $\rho_3 \notin H$, since $H$ has order $2k$ and $k$ is odd. We conclude that $(\psl_2(q), (\rho_0,\rho_1,\rho_2,\rho_3))$ is a rank $4$ extension of $\M$. 
\end{proof}

Our next task is to construct a infinite family of reflexible $4$-maniplexes contained in each of the classes 1 and $2$. We will make use of Dirichlet's Theorem stating that if $a$ and $b$ are two coprime integers, then there are infinitely many prime numbers congruent to $a$ modulo $b$.

\subsection{Class 1}

\begin{lemma}
\label{lem:class1}
Let $p$ be a prime such that $p \equiv 1$ $(\mod 12)$. Then there exists a reflexible map $\M$ of type $\{3,p\}$ with automorphism group $\psl_2(p)$ admitting a reflexible non-orientable extension $\bar{\M}$ such that $\Aut(\bar{\M}) = \psl_2(p)$. Moreover, $\bar{\M}$ belongs to Class $1$.
\end{lemma}

\begin{proof}
We will show that $\M$ exists by providing an explicit set of generators for its automorphism group. Let $a \in \text{GF}(p)$ be such that $a^2 = -1$ (since $p \equiv 1$ $(\mod 4)$) such an element always exists) and consider the following matrices in $\textrm{SL}_2(p)$.

\begin{align*}
\overline{\rho_0} =
\begin{pmatrix}
0 & 1 \\
-1 & 0\\
\end{pmatrix}
\qquad
\overline{\rho_1} =
\begin{pmatrix}
a & 0 \\
1 & -a\\
\end{pmatrix}
\qquad
\overline{\rho_2} =
\begin{pmatrix}
-a & 0 \\
0 & a\\
\end{pmatrix}
\qquad
\end{align*}

For each $\overline{\rho_i}$, let $\rho_i$ be the corresponding image in $\psl_2(p)$ under the quotient projection $\pi\colon \textrm{SL}_2(q) \to \psl_2(q)$. Observe that $(\rho_i)^2 = 1$ for all $i \in \{0,1,2\}$ and that $(\rho_0\rho_2)^2 = (\rho_0\rho_1)^3 = 1$. Moreover, for an integer $k$ we have

\begin{align*}
(\overline{\rho_1}\overline{\rho_2})^k =
\begin{pmatrix}
1 & 0 \\
-ka & 1\\
\end{pmatrix}.
\qquad
\end{align*}

and thus $\rho_1\rho_2$ has order $p$. We must show that group $G := \langle \rho_0, \rho_1, \rho_2 \rangle$ is isomorphic to $\psl_2(p)$. Since $G \leq \psl_2(q)$ and $G$ is non-abelian, it suffices to show that $G$ is not of any of the four possible types of proper non-abelian subgroups of $\psl_q(p)$ given in Lemma \ref{lem:dick}.  

First, note that since $p$ is prime, $\psl_2(p)$ has no proper projective linear subgroups. Now, $G$ contains a dihedral subgroup $\langle \rho_0, \rho_1 \rangle$ of order $6$, and thus $G$ is not a semidirect product of an elementary abelian group with a cyclic group. To see that $G$ cannot be dihedral, recall that in a dihedral group the product of two non-central involutions that commute must be the central element of the group. However, $\rho_0\rho_2$ is not central in $G$ since it does not commute with $\rho_1$. Finally, $G$ is not isomorphic to $\rA_5$ or $\rS_4$ since $\rho_1\rho_2$ has order $p \geq 13$. Then $G = \psl_2(p)$ and $\M := (G,(\rho_0,\rho_1,\rho_2))$ is a reflexible map of type $\{3,p\}$ with automorphism group $\psl_2(p)$. 
 
Now, $3$ divides $(p-1)/2$ and by Theorem \ref{the:ext} there exists a $\rho_3 \in G$ such that $(G,(\rho_0,\rho_1,\rho_2,\rho_3))$ is a reflexible $4$-maniplex with automorphism group $G$. It remains to show that $(G,(\rho_0,\rho_1,\rho_2,\rho_3))$ belongs to Class 1. Clearly $G_0 = \langle \rho_1, \rho_2, \rho_3 \rangle$ is not isomorphic to $\rA_5$ or $\rS_4$ since it contains an element of order $p$ and by Theorem \ref{the:facet}, $G_0 \cong \psl_2(p)$ and thus $G_0 = G_3 = G$, as desired.
\end{proof}

Since there are infinitely many prime numbers congruent to $1$ modulo $12$, we obtain the following corollary.

\begin{corollary}
There are infinitely many reflexible $4$-maniplexes in Class 1.
\end{corollary}

\subsection{Class 2}

\begin{proposition}
\label{prop:class2}
Let $p \neq 19$ be a prime and suppose that for $\alpha \in \{-1,1\}$, $p \equiv \alpha$ $(\mod 20)$. Then there exists a reflexible $4$-maniplex $\M$ with automorphism group $\psl_2(p)$. 
\end{proposition}

\begin{proof}
Let $G = \psl_2(p)$ and $n = (p-\alpha)/2$. Observe that, because of our choice of $p$, $G$ admits a maximal subgroup $A$ isomorphic to $\textrm{A}_5$. It is well known that $A$ can be generated by three involutions $\rho_0$, $\rho_1$ and $\rho_2$ such that $\rho_0\rho_2 = \rho_2\rho_0$ and the orders of $\rho_0\rho_1$ and $\rho_1\rho_2$ are $5$ and $3$ respectively \cite{cherk}. The group $H = \langle \rho_0, \rho_1 \rangle$ is isomorphic to $\rD_{10}$ and is thus contained in a maximal subgroup $D$ of $G$ isomorphic to $\rD_{2n}$. Since $n$ is even, $D$ contains a non-trivial central involution $\rho_3$. Now, the centralizers of two non-commuting involutions in $\rA_5$ must have trivial intersection and thus $\rho_3 \notin A$. Then, by the maximality of $A$, we see that $G = \langle \rho_0, \rho_1, \rho_2, \rho_3 \rangle$ and thus $(G,(\rho_0,\rho_1,\rho_2,\rho_3))$ is a reflexible $4$-maniplex with automorphism group $G = \psl_2(p)$. 
\end{proof}

The reflexive $4$-maniplex constructed in Proposition \ref{prop:class2} belongs to Class 2. To see this, first note that $G_0 := \langle \rho_1, \rho_2, \rho_3 \rangle$ must be $\psl_2(p)$. Indeed, by Theorem \ref{the:facet}, $G_0$ is isomorphic to $\rS_4$, $\rA_5$ or $\psl_2(p)$. Since $p \notin \{11,19\}$ and $A$ is alternating, by Proposition \ref{pro:class3}, $G_0$ cannot be isomorphic to $\rS_4$ or $\rA_5$. Hence, $G_0 = \psl_q(p)$ and by construction $G_3 := \langle \rho_0, \rho_1, \rho_2 \rangle$ is isomoprhic to $\rA_5$. Therefore, $\M$ is in Class $2$. Moreover, $\M$ is an extension of the reflexible map $\M_3 := (G_3,(\rho_0, \rho_1, \rho_2))$ of type $\{5,3\}$ isomorphic to the hemi-dodecahdron. Since there are infinitely many primes congruent to $\pm 1$ modulo $20$, the construction of Proposition \ref{prop:class2} gives us infinitely many maniplexes

\begin{corollary}
There are infinitely many reflexible $4$-maniplexes in Class 2.
\end{corollary}


\section{Remarks on Suzuki groups and $\pgl_2(q)$}\label{sec:remarks}

Some of the strategies used in this paper may also prove useful when looking at other families of groups. Consider, for instance, the case of Suzuki groups (see \cite{Tits1961} for definition and details). For an odd power of 2 $q = 2^{2n+1}$, the Suzuki group $\Sz(q)$ can be seen as a doubly-transitive permutation group acting on a set $\Omega$ of size $q^2+1$. With this action, all involutions of $\Sz(q)$ fix exactly one point of $\Omega$ and two involutions commute if and only if they have same fixed point. Hence, the argument used in Lemma \ref{lem:qeven} shows that there are no reflexible maniplexes of rank $n$ with full automorphism group isomorphic to a Suzuki group if $n > 3$.

Although left out of the scope of this paper, maniplexes having a projective general linear group as their full automorphism group can also be tackled by following the same strategy used here. Indeed, since for every prime power $q$ we have $\psl_2(q) \leq \pgl_2(q) < \psl_2(q^2)$, the subgroup structure of $\pgl_2(q)$ is well understood (see for instance \cite{pglpoly} for a detailed list). In particular a proper non-abelian subgroup of $\pgl_2(q)$ must be either $\rA_4$, $\rS_4$, $\rA_5$ or $\psl_2(q)$. Thus, up to duality, the facet group of a maniplex of rank $n$ with full automorphism group $\pgl_2(q)$ must be isomorphic to $\rS_4$, $\rA_5$, $\psl_2(q)$ or $\pgl_2(q)$. As $\pgl_2(q)$ does not admit a rank $n$ string group representation for $n > 4$, the maximum rank for a reflexible maniplex with  automorphism group $\pgl_2(q)$ is $4$. A computer assisted search for rank $4$ maniplexes with automorphism group $\pgl_2(q)$ for small values of $q$ suggests that the $4$-simplex is the only $4$-maniplex having proper subgroups of $\psl_2(q)$ as its vertex and facet groups (it is known that it is the only regular abstract polytope of rank $4$ with automorphism group $\pgl_2(q)$ \cite{pglpoly}). For all other $4$-maniplexes found (up to duality), the facet group is isomorphic to $\pgl_2(q)$ while the vertex group is one of $\rS_4$, $\rA_5$, $\psl_2(q)$ or $\pgl_2(q)$, with all four possibilities occurring. We conjecture that there is an infinite family of rank $4$ maniplexes for all four of the aforementioned possibilities for the vertex group.


\begin{thebibliography}{SK}

\bibitem{cherk} M. Cherkassoff, D. Sjerve, On groups generated by three involutions, two of which commute, {\em The Hilton symposium 1993} (Montreal, PQ), CRM Proc. Lecture Notes, vol. 6, Amer. Math. Soc., Providence, RI, 1994, pp 169-185. MR1290589

\bibitem{Conder2008}
M. Conder, P. Poto{\v{c}}nik, and J. {\v{S}}ir{\'a}{\v{n}}.
\newblock Regular hypermaps over projective linear groups.
\newblock {\em J. Aust. Math. Soc.}, 85(2):155--175, 2008.

\bibitem{Cox1982} H.S.M.~Coxeter. 
\newblock Ten toroids and fifty-seven hemi-dodecahedra.
\newblock {\em Geom.\ Dedicata,\/} 13:87--99, 1982.

\bibitem{stg} G. Cunningham, M. Del Río-Francos, I. Hubard, M. Toledo, Symmetry type graphs of polytopes
and maniplexes, {\em Ann. Comb.} {\bf 19} (2015), p. 243--268.

\bibitem{dickson} L.E. Dickson, Linear groups: with an exposition of the Galois field theory, Dover publications Inc., New York, 1958.

\bibitem{mani2} I. Douglas, I. Hubard, D. Pellicer, S. Wilson, The twist operator on maniplexes, in: Springer Contributed Volume on Discrete Geometry and Symmetry, in: Springer Proc. Math. Stat., vol. 234, 2018, pp. 127–145.

\bibitem{polymani} J. Garza-Vargas, I. Hubard, Polytopality of maniplexes, {\em Discrete Math.} 341:2068--2079, 2018.

\bibitem{michael} M. Giudici, Maximal subgroups of almost simple group with socle $\psl(2,q)$, 	arXiv:math/0703685.

\bibitem{Gr1977} B.~Gr\"unbaum.
\newblock Regularity of graphs, complexes and designs.
\newblock {\em Probl\`{e}mes Combinatoires et Th\'{e}orie des Graphes}, volume 260 of
{\em Coll.\  Int.\ C.N.R.S.\/}, pages 191--197, Orsay, 1977. 

\bibitem{Hup67}
B.~Huppert.
\newblock {\em Endliche {G}ruppen, {I}}.
\newblock Die Grundlehren der Mathematischen Wissenschaften, Band 134. Springer-Verlag, Berlin, 
1967.


\bibitem{LV2005}
D.~Leemans and L.~Vauthier.
\newblock An atlas of abstract regular polytopes for small groups.
\newblock {\em Aequationes Math.}, 72:313--320, 2006.

\bibitem{pslpoly} D. Leemans, E. Schulte, Groups on type $L_2(q)$ acting on polytopes, {\em Adv. Geom.} 7:529--539, 2007.

\bibitem{pglpoly} D. Leemans, E. Schulte, Polytopes with groups of type $\pgl_2(q)$, {\em Ars Math. Contemp.} {\bf 2}, (2009), pp. 163--171.

\bibitem{LV2005}
D.~Leemans and L.~Vauthier.
\newblock An atlas of abstract regular polytopes for small groups.
\newblock {\em Aequationes Math.,\/} 72:313-320, 2006.

\bibitem{combimaps} G. A. Jones and J. S. Thornton, Operations on maps and outer automorphisms, {\em 0J. Combin. Theory
Ser. B} 35:93--103, 1983.

\bibitem{macbeath}
A. M. Macbeath, Generators of the linear fractional groups, {\em Number Theory, Proc. Sympos. Pure Math., Houston, TX,} Vol. \textrm{XII}, (American Mathematical Society, Providence, RI, 1967), p. 14--32.

\bibitem{2orbit} D. Pellicer, P. Poto\v{c}nik, M. Toledo, An existence result on two-orbit maniplexes, {\em J. Combin. Theory Ser. A}, {\em  166} (2019), 226--253. 

\bibitem{dansteve} D. Pellicer, S. Wilson, Rotary one-facet maniplexes, {\em Art Discret. Appl. Math.} {\bf 4} (2021), P3.02.

\bibitem{sah}
C.-H. Sah, Groups related to compact Riemann surfaces, {\em Acta Math.}, 123:13--42, 1969.

\bibitem{Tits1961} J.~Tits.
\newblock Les groupes simples de Suzuki et de Ree.  
\newblock {\em S\'eminaire N. Bourbaki}, exp. 210, p. 65--82, 1960--1961.

\bibitem{Van85} P.~Vanden Cruyce.
\newblock G\'eom\'etries des groupes PSL(2,q).  
\newblock {\em PhD. Thesis, Universit\'e Libre de Bruxelles}, 1985.

\bibitem{manisteve} S. Wilson, Maniplexes Part I: maps, polytopes, symmetry and operators, {\em Symmetry} 4:265--275, 2012.



\end{thebibliography}
\end{document}